\documentclass[11pt]{article}
\usepackage[T1]{fontenc}
\usepackage[a4paper]{anysize}\marginsize{2.0cm}{2.0cm}{1.0cm}{1.0cm}
\pdfpagewidth=\paperwidth \pdfpageheight=\paperheight
\usepackage{pgf,tikz,hyperref}
\usetikzlibrary{arrows}
\usepackage{amsfonts,amssymb,amsthm,amsmath,eucal}
\usepackage{graphicx}
\usepackage{ltablex}
\usepackage{caption, booktabs}
\usepackage{longtable}

\newcommand{\Loc}{\operatorname{Loc}}
\newcommand{\coef}{\operatorname{coef}}

\theoremstyle{plain}
\newtheorem{thm}{Theorem}[section]
\newtheorem{theorem}[thm]{Theorem}

\newtheorem{lemma}[thm]{Lemma}
\newtheorem{proposition}[thm]{Proposition}

\newtheorem{custom}{Theorem}

\newtheorem{subtheorem}{Theorem}
\newenvironment{theorem2}[1]{%
  \renewcommand\thesubtheorem{#1}%
  \subtheorem
}{\endsubtheorem}

\theoremstyle{definition}

\newtheorem{thevarthm}[thm]{\varthmname}

\newenvironment{varthm*}[1]{\trivlist\item[]{\bf #1.}\it}{\endtrivlist}

\newcommand\be{\begin{eqnarray*}}
\newcommand\ee{\end{eqnarray*}}

\renewcommand\P{\mathbb P}

\newcommand\newop[2]{\def#1{\mathop{\rm #2}\nolimits}}
\newop\edim{edim}
\newop\Zeroes{Zeroes}
\newop\Jac{Jac}
\newop\Ass{Ass}
\newop\SL{SL}
\newop\PGL{{\P}GL}
\newop\Km{Km}
\newop\reg{reg}

\newcolumntype{L}{>{$}l<{$}}

\newcommand\keywords[1]{{\renewcommand\thefootnote{}\footnotetext{\textit{Keywords:} #1.}}}
\newcommand\subclass[1]{{\renewcommand\thefootnote{}\footnotetext{\textit{Mathematics Subject Classification (2010):} #1.}}}

\begin{document}

\author{\L{}ukasz Merta \and Maciej Zi\k{e}ba}
\title{Sextactic and type-9 points on the Fermat cubic\\ and associated objects}
\date{\today}
\maketitle
\thispagestyle{empty}

\begin{abstract}
In this note we study line and conic arrangements associated to sextactic and type $9$ points on the Fermat cubic $F$ and we provide explicit coordinates for each of the $72$ type $9$ points on $F$.
\end{abstract}
\keywords{Fermat cubic, torsion points, line arrangement, conic arrangement}
\subclass{MSC 14C20 \and MSC 14N20 \and MSC 13A15}

%*****************************************************************************
\section{Introduction}
Our research is motivated by recent results of Szemberg and Szpond \cite{SzeSzp24} on sextactic points on the Fermat cubic. Recall that for an arbitrary smooth plane curve $C$ of degree at least $3$ a point $P\in C$ is called a \emph{sextactic point}, if there exists an irreducible hyperosculating conic $D$ at $P$. This means that the order of tangency  between $C$ and $D$ at $P$ is at least $6$. Note, that it was proved by Cayley already in 1859 that for every point $P\in C$ there exists a conic with tangency order at least $5$ at $P$, see \cite{Cay1859}. Thus sextactic points are a natural generalization of flex points of a plane curve $C$, i.e., points at which the order of tangency of the tangent line is strictly bigger than $2$.

This notion has been further generalized to type $9$ points by Gattazzo \cite{Gat79}. These are points $P$ on $C$ such that there exists an irreducible cubic curve, which intersects $C$ at $P$ with multiplicity $9$. 

If $C$ is a curve of degree $3$, then the flexes are the only common points of $C$ and the line tangent at the flex. Similarly, the sextactic points are the only intersection points of the cubic with the hyperosculating conic. Note that the flex points are not sextactic because for them the hyperosculating conic is just the tangent line taken as a double line, so that it is not irreducible. The same remains true for type $9$ points.

In \cite[Problem 5.3]{SzeSzp24} Szemberg and Szpond asked for explicit coordinates of the type $9$ points on the Fermat cubic. A solution to this problem is our first result (we provide explicit formulas in Section \ref{sec: type 9}).
\begin{custom} ({see Section \ref{sec: type 9}} for details) \label{thm: coordinates}
Let $F$ be the Fermat cubic $x^3+y^3+z^3=0$. Then the type $9$ points on $F$ have coordinates defined in the field $\mathbb{Q}(\alpha,\beta)$, where $\alpha=\sqrt[3]{3}$ and $\beta$ is a primitive root of unity of order $9$.
\end{custom}
There are $72$ type $9$ points on $F$ (a result already known to Gattazzo) and their union is indeed a complete intersection of type $(3,24)$ as stipulated in \cite[Problem 5.2]{SzeSzp24}.
Our two other results concern geometry of certain line arrangements associated naturally to sextactic and type $9$ points.
\begin{custom} ({see Theorem \ref{thm: tangent_s}})\label{thm: tangents at sex}
%Let $S_1,\ldots,S_{27}$ be the sextactic points on the Fermat cubic $F$ and let $L_i$ be the tangent line at $S_i$ for $i=1,\ldots,27$. Then the residual intersection point of $L_i$ with $F$ is one of its flex points.
For any sextactic point $S \in F$, the residual intersection point of the line tangent to $F$ at $S$ is one of its flex points.
\end{custom}

The above observation provides additionally a partition of the set of tangent lines $\left\{L_1,\ldots,L_{27}\right\}$ in triads, i.e., there are exactly three of these lines passing through any of the flex points.

\begin{custom} ({see Theorem \ref{thm: tangent_t9}})
For any type $9$ point $T\in F$, the residual intersection point of the tangent line to $F$ at $T$ is another point of type $9$.
\end{custom}

We also consider conics which intersect $F$ in two points (which are either two sextactic points or two points of type $9$), so that the intersection multiplicity between $F$ and the conic at each of these points is exactly $3$.

\begin{custom} ({see Section \ref{sec: conics} for details})

i) For any sextactic point $S\in F$, there are exactly $8$ conics which intersect $F$ with multiplicity $3$ at $S$ and another sextactic point $S'$. All these $108$ conics intersect with multiplicity either $6$ or $9$ in $54$ ''shadow'' points off the Fermat cubic and in $486$ ordinary double points.

ii) For any type $9$ point $T\in F$, there are exactly $9$ conics which intersect $F$ with multiplicity $3$ at $T$ and another type $9$ point $T'$. All these $324$ conics intersect with multiplicity $9$ in $72$ ''shadow'' points off the Fermat cubic and in $1944$ ordinary double points.
\end{custom}
We illustrate how these ''shadow'' points mentioned above are distributed and we also give their explicit coordinates over respective fields.

\section{Type $9$ points on the Fermat cubic}\label{sec: type 9}

Type 9 points on a Fermat cubic (as well as any other cubic) are precisely the $9$-torsion points on $F$, which are not $3$-torsion points. This is a consequence of Abel's Theorem, which we state in a simpler form below:

\begin{theorem}[Abel's Theorem]
    Let $E$ be a smooth complex elliptic curve with distinguished point $0$ embedded in the complex projective plane. Then a divisor $D = \sum d_iP_i$ on $E$ is a (scheme theoretic) intersection of $E$ with another plane curve $C$ of degree $c$ if and only if $\sum d_i = 3c$ and $\sum d_iP_i = 0$ in the group law on $E$.
\end{theorem}

Therefore, to calculate the explicit coordinates of type $9$ points on the Fermat cubic, we need to calculate its $9$-torsion points. This can be done by using the division polynomials. Let us recall that for a given elliptic curve $y^2 = x^3 + Ax + B$ such polynomials $(\psi_n)_{n \in \mathbb{N}}$ are defined recursively in the following way:
\begin{align*}
\psi_0 &= 0, \\
\psi_1 &= 1, \\
\psi_2 &= 2y, \\
\psi_3 &= 3x^4 + 6Ax^2 + 12Bx - A^2, \\
\psi_4 &= 4y(x^6 + 5Ax^4 + 20Bx^3 - 5A^2x^2 - 4ABx - 8B^2 - A^3), \\
\psi_{2k+1} &= \psi_{k+2}\psi_k^3 - \psi_{k-1}\psi_{k+1}^3 \; \text{for } k \geq 2,\\
\psi_{2k} &= \frac{\psi_k}{2y} (\psi_{k+2}\psi_{k-1}^2 - \psi_{k-2}\psi_{k+1}^2) \; \text{for } k > 2.
\end{align*}
The roots of $\psi_9$ (which becomes a polynomial in one variable $x$ after substitution $y^2 = x^3 + Ax + B$) are precisely the $x$-coordinates of the $9$-torsion points of that elliptic curve (except of the origin point). In order to use this method for Fermat cubic, we need to write it in the Weierstrass form first. This can be done under the following substitutions
\begin{equation}\label{subst}
\tilde{x} = \frac{-12z}{x+y}, \;\; \tilde{y} = \frac{36(x-y)}{x+y},
\end{equation}
which result in the equation
\[\tilde{y}^2 = \tilde{x}^3 - 432.\]
We therefore have $A = 0$ and $B = -432$. The next step is to find $\gamma \in \mathbb{C}$ such that the obtained polynomial $\psi_9$ of degree $40$ splits into linear factors over the field $\mathbb{Q}(\gamma)$. It turns out that we can use $\gamma = \alpha + \beta$, where $\alpha = \sqrt[3]{3}$ and $\beta$ is the primitive root of unity of order $9$. This is an algebraic number of degree $18$, satisfying the equality
\[\gamma^{18} - 15\gamma^{15} + 177\gamma^{12} - 578\gamma^9 + 6747\gamma^6 + 642\gamma^3 + 343 = 0.\]
Therefore, the exact coordinates of all type $9$ points on the Fermat cubic can be defined in the field $\mathbb{Q}(\gamma) = \mathbb{Q}(\alpha, \beta)$. By factoring the polynomial $\psi_9$ over this field, we obtained $40$ $x$-coordinates, hence $80$ different points (there are $2$ values of $y$ for each $x$). Then we applied the inverse of \eqref{subst}, which is
\begin{equation}\label{subst_inv}
x = 36 + \tilde{y}, \;\; y = 36 - \tilde{y}, \;\; z = -6x,
\end{equation}
to obtain the list of coordinates of $80$ points on projective plane. Naturally, that list still contained eight $3$-torsion points of $F$, which we had to remove -- these were precisely the points lying on the Hessian curve of $F$, namely $xyz = 0$. That way we finally obtained the list of $72$ points.

The explicit coordinates of the first $12$ points of type $9$ on the Fermat cubic are written below. By permuting the coordinates of these $12$ points, we can obtain the coordinates of the remaining $60$ points from our list.

% Najprościej zapisane współrzędne jakie tylko udało mi się uzyskać.
% W 9 punktach postanowiłem dać jako pierwszą współrzędną 3, a nie 1, bo w pozostałych współrzędnych wychodziły ułamki 1/3 i 2/3, teraz wygląda to po prostu lepiej i czytelniej.
% Skoro pozostałe punkty to permutacje poniższych, nie ma sensu chyba wypisywać ich wszystkich?
\vspace{-10pt}
\begin{align*}
T_1 &= [1 : \beta : \beta^2], \\
T_2 &= [1 : \beta^2 : \beta^4], \\
T_3 &= [1 : \beta : \beta^5], \\
T_4 &= [3 : \alpha\beta(\beta^4 + 2 \beta^3 - \beta + 1) : -\alpha^2(\beta^2 + \beta + 1)(\beta^3 - \beta + 1)], \\
T_5 &= [3 : -\alpha\beta(2\beta^4 + \beta^3 + \beta + 2) : \alpha^2(\beta+1)(\beta-1)^2], \\
T_6 &= [3 : \alpha\beta(\beta^4 - \beta^3 - \beta - 2) : -\alpha^2\beta^2(\beta^2 + \beta + 1)], \\
T_7 &= [3 : \alpha\beta(\beta - 1)(\beta^3 + 2) : -\alpha^2\beta(\beta^2 + \beta + 1)(\beta^2 - \beta + 1)], \\
T_8 &= [3 : \alpha^2(\beta - 1)(\beta + 1)(\beta^3 + \beta + 1) : \alpha\beta(\beta^4 - \beta^3 - \beta - 2)], \\
T_9 &= [3 : \alpha\beta(\beta-1)^2(\beta^2 + \beta + 1) : -\alpha^2\beta(\beta^2 + \beta + 1)(\beta^2 - \beta + 1)], \\
T_{10} &= [3 : \alpha\beta(\beta - 1)^2(\beta^2 + \beta + 1) : alpha^2(\beta - 1)(\beta^3 + \beta^2 + 1)], \\
T_{11} &= [3 : \alpha\beta(\beta^4 + 2\beta^3 + 2\beta + 1) : -\alpha^2\beta^2(\beta^2 + \beta + 1)], \\
T_{12} &= [3 : -\alpha^2(\beta^2 + \beta + 1)(\beta^3 - \beta + 1) : -\alpha\beta(2\beta^4 + \beta^3 + \beta + 2)].
\end{align*}

Since flex points on $F$ are a complete intersection of $F$ and its hessian ($H(F) = xyz$), while sextactic points on $F$ are a complete intersection of $F$ and its second hessian (see \cite{SzeSzp24} for details), given by 
$$H_2(F) = (x^3 - y^3)(y^3 - z^3)(x^3 - z^3),$$
it is natural to expect that the set of points of type $9$ on $F$ is a complete intersection of $F$ and some curve of degree $24$. With the explicit coordinates of type $9$ points, it is now possible to compute the standard basis of an ideal $I_T = I(T_1, T_2, \dots, T_{72}) \subset \mathbb{Q}(\alpha, \beta)[x,y,z]$ of those points. We obtained
$$I_T = (F, Q(x,y)) = (F, Q(x,z)) = (F, Q(y,z)),$$
where
\begin{align*}
Q(x,y) &= x^{24} + 4x^{21}y^3-17x^{18}y^6-65x^{15}y^9-89x^{12}y^{12}-65x^9y^{15}-17x^6y^{18}+4x^3y^{21}+y^{24} \\
&= (x^9 - 3x^6y^3 - 6x^3y^6 - y^9)(x^6 + x^3y^3 + y^6)(x^9 + 6x^6y^3 + 3x^3y^6 - y^9).
\end{align*}
Moreover, it turns out that $Q(x,y)$ splits into $24$ linear factors over $\mathbb{Q}(\alpha, \beta)$.
% TODO: Add equetion of curve of degree 24, add comment about factorazing it to 24 lines

% TODO: Check double conics for sextatic points. 

\section{Tangents to the Fermat cubic}

In this section, we investigate the arrangements of lines which are tangent to the Fermat cubic $F$ at sextactic points and type $9$ points.

%For each sextactic point $S_i$ on $F$, let $L_i$ be the tangent line at $S_i = [a_i:b_i:c_i]$. Considering the Fermat cubic $F$ and its partial derivatives, we can compute the equations of tangent lines at each sextactic point. These equations will involve the coordinates of the sextactic points and will be specific to each $L_i$. The equation of a line passing through a point $S_i$ can be expressed as:
%\[
%L_i: \quad a_i^2 \cdot x + b_i^2 \cdot y + c_i^2 \cdot z = 0.
%\]

For each point $P = [a_i : b_i : c_i]$ on $F$, the tangent line $L$ is given by equation
\[
L_i: \quad a_i^2 \cdot x + b_i^2 \cdot y + c_i^2 \cdot z = 0.
\]

%Using the equations of the tangent lines $L_i$, we can find the points of intersection between these lines and the Fermat cubic $F$. The saturation process allows us to determine these intersection points.

If $P$ is a flex point on $F$, then the tangent line has multiplicity $3$, hence it cannot intersect with $F$ in any other point besides $P$. If it is not a flex point, the tangent line at $P$ intersects Fermat cubic at exactly one additional point besides $P$. Let us recall here that the flex points on Fermat cubic are precisely the points lying on lines $x = 0$, $y = 0$ and $z = 0$, with coordinates
\begin{gather*}
    [1 : -1 : 0],\ [1 : 0 : -1],\ [0 : 1 : -1], \\
    [1 : -\varepsilon : 0],\ [1 : 0 : -\varepsilon],\ [0 : 1 : -\varepsilon], \\
    [1 : -\varepsilon^2 : 0],\ [1 : 0 : -\varepsilon^2],\ [0 : 1 : -\varepsilon^2],
\end{gather*}
where $\varepsilon$ is the primitive root of unity of order $3$. 

Since sextactic points and type $9$ points are not flex points, we can consider their respective residual intersection points. Naturally, for a line tangent to $F$ at a sextactic point and a type $9$ point, we expect the residual intersection point to be a $6$-torsion point and a $9$-torsion point, respectively -- this is the consequence of the group law on $F$. 

We calculated the coordinates of these intersection points and we obtained the following results.

\begin{theorem} \label{thm: tangent_s}
    Each line tangent to Fermat cubic at a sextactic point intersects the cubic in a flex point. Moreover, there are three such lines going through each flex point.
\end{theorem}

\begin{theorem2}{\ref{thm: tangent_s}'}{\label{thm: tangent_t9}}
    Each line tangent to Fermat cubic at a point of type $9$ intersects the cubic in another point of type $9$.
\end{theorem2}

We can study the arrangement of lines tangent to Fermat cubic at points of type $9$ in more detail. Denote by $\mathcal{T}$ the set of points of type $9$ on Fermat cubic. We consider a mapping $\nu : \mathcal{T} \mapsto \mathcal{T}$, such that for each $T \in \mathcal{T}$ the line tangent to Fermat cubic at $T$ intersects the cubic in another point $\nu(T)$. For each point $T \in \mathcal{T}$ we can then consider a sequence of points
$$T,\ \nu(T),\ \nu(\nu(T)),\ \dots$$
and see how long it takes for this sequence to repeat itself. It turns out that we have
$$\nu(\nu(\nu(T))) = T$$
for all $T \in \mathcal{T}$. Therefore the set $\mathcal{T}$ can be divided into $24$ disjoint subsets of the form $\{T,\ \nu(T),\, \nu(\nu(T))\}$ for some $T \in \mathcal{T}$.

% By using the equations of tangent lines $L_i$, we explored where these lines intersect the Fermat cubic $F$. Interestingly, when we looked at the points where each tangent line crosses the Fermat cubic for the third time, we found only nine such points. Each of these points is shared by three tangent lines, and we call them $P_i^{(j)}$, where $i$ represents the sextactic point, and $j$ indicates one of the three tangent lines.
% What's fascinating is that each of these points not only lies on the Fermat cubic but also aligns with the Hessian of the Fermat cubic. This alignment reveals that these points are flex points on the Fermat curve.

\section{Arrangements of conics}\label{sec: conics}

In this section, we consider conics $C$, passing through two distinct points on Fermat cubic $F$, such that the intersection multiplicities of $C$ and $F$ at these points are both equal to $3$. Our aim is to find all such conics passing through either two sextactic points or two points of type $9$ on $F$, then study their intersection patterns.

Before we start, note that from Abel's theorem, we can obtain a simple criterion for such conics to exist for two given points on an elliptic curve:

\begin{proposition} \label{prop: thirdpoint}
    Let $P_1$ and $P_2$ be distinct points on a smooth complex elliptic curve $E$. Then there exists a conic $C$  tangent to $E$ at those two points with multiplicity $3$ in each point if and only if $P_1 + P_2$ is a $3$-torsion point on $E$.
\end{proposition}

Indeed, an intersection of $E$ with a conic $C$ is a divisor $D = 3P_1 + 3P_2$ if and only if $3(P_1 + P_2) = 0$ in the group law on $E$. The latter condition means that the line passing through $P_1$ and $P_2$ must intersect $E$ in the third point, which is a $3$-torsion point (and also a flex). Therefore, since an elliptic curve  $E$ has $9$ flex points, we can expect $9$ such conics passing through a generic point $P$ on $E$.

We are now going to briefly describe the algorithm we used to find such conics. First, we introduce a lemma which we are going to use later. Here, $I_P(f,g)$ denotes the intersection multiplicity of affine curves $f=0$ and $g=0$ at point $P$.

\begin{lemma}\label{lem: intersect}
    Let $f, g \in \mathbb{K}[x,y]$ be two affine curves such that $f(0,0) = g(0,0) = 0$. If $f$ and $g$ can be written as
    \begin{align*}
        f &= ax + \text{(remaining terms)}, \\
        g &= by^k + \text{(higher degree terms)}
    \end{align*}
    with $a, b \in \mathbb{K}$ and $k \in \mathbb{N}$, then $I_{(0,0)}(f,g) \leq k$.
\end{lemma}

Before we show the proof of this lemma, we introduce some notation and some additional lemmas, based on \cite{Gra94}. Let $f \in \mathbb{K}[x]$, where $x = (x_1, x_2, \dots, x_n)$. Let $<$ be any monomial order, i.e.\ a total order on the set of monomials, respecting multiplication. We denote by $L(f)$ the leading monomial of $f$ with respect to the given monomial order, without the coefficient.  Let $\Loc_<\mathbb{K}[x]$ be a localization of $\mathbb{K}[x]$ with respect to the multiplicative set
\[S_< = \{ 1 + g : g = 0 \text{ or } g \in \mathbb{K}[x]\setminus\{0\} \text{ and } 1 > L(g)\}.\]
For an ideal $J \subset \mathbb{K}[x]$, let $L(J)$ be the ideal generated by $\{L(f) : f \in J\}.$ We say that the set $\{f_1, f_2, \dots, f_s\} \subset J$ is a \emph{standard basis} of $J$ if $\{L(f_1), L(f_2), \dots, L(f_s)\}$ generates $L(J)$.

We will use the following two results from \cite{Gra94}.

\begin{proposition}[{\cite[Proposition 1.4]{Gra94}}]\label{proposition: std} If $\{f_1, f_2, \dots, f_s\}$ is a standard basis of $J$, then 
\[J \Loc_<\mathbb{K}[x] = (f_1, f_2, \dots, f_s)\Loc_<\mathbb{K}[x].\]
\end{proposition}

\begin{proposition}[based on {\cite[Corollary 5.4]{Gra94}}]\label{proposition: basis} Let either $<$ be a well-order or $\Loc_<\mathbb{K}[x]/J\Loc_<\mathbb{K}[x]$ a finite-dimensional vector space over $\mathbb{K}$. Then the monomials in $\mathbb{K}[x]\setminus L(J)$ represent a $\mathbb{K}$-basis of $\Loc_<\mathbb{K}[x]/J\Loc_<\mathbb{K}[x]$.
\end{proposition}

\begin{proof}[Proof of Lemma \ref{lem: intersect}]
    Let $J = (f,g)$. Consider the following order of monomials:
    \[1 > x > y > x^2 > xy > y^2 > \dots\]
    Note that with this ordering, the set $\Loc_<\mathbb{K}[x,y]$ is actually equal to $\mathbb{K}[x,y]_{(x,y)}$, i.e.\ a localization with respect to the maximal ideal. Therefore
    \[I_{(0,0)}(f,g) = \dim \mathbb{K}[x,y]_{(x,y)}/J\mathbb{K}[x,y]_{(x,y)} = \dim \Loc_<\mathbb{K}[x,y]/J\Loc_<\mathbb{K}[x,y].\]
    We can choose additional forms $h_1, \dots, h_s$ so that $\{f, g, h_1, \dots, h_s\}$ is a standard basis of $J$. By Proposition \ref{proposition: std} we know that $\{f, g, h_1, \dots, h_s\}$ generates $J$ in $\Loc_<\mathbb{K}[x,y]$. By Proposition \ref{proposition: basis} we have that the $\mathbb{K}$-basis of $\Loc_<\mathbb{K}[x,y]/J\Loc_<\mathbb{K}[x,y]$ is contained in the set of monomials in $\mathbb{K}[x,y] \setminus L(J)$. Since $\{f, g, h_1, \dots, h_s\}$ is a standard basis, the ideal $L(J)$ is generated by $\{L(f), L(g), L(h_1), \dots, L(h_s)\}$. Furthermore, since $L(f) = x$ and $L(g) = y^k$, the set of monomials in $\mathbb{K}[x,y] \setminus L(J)$ is contained in
    \[\{1, y, y^2, \dots, y^{k-1}\},\]
    hence $\dim \Loc_<\mathbb{K}[x,y]/J\Loc_<\mathbb{K}[x,y] \leq k$, which completes the proof.
\end{proof}

In order to find conics described at the beginning of this section, we use the equation of Fermat cubic $F = x^3 + y^3 + z^3 = 0$ and a general cubic equation $C = ax^2 + by^2 + cz^2 + dxy + exz + fyz = 0$. In order to use Lemma \ref{lem: intersect}, we need to switch to the affine case, so we assume $x = 1$ in both equations. Let $P = [1 : y' : z']$ be a point on $F$, which is either a sextactic point or a point of type $9$. Assume that $C$ intersects $F$ at $P$ with multiplicity $3$. First, we use substitutions
\[y \leftarrow y + y', \quad z \leftarrow z + z'\]
to shift both $C$ and $F$ to the point $(0,0)$. Since $P$ lies on $C$, we have
\begin{equation}\label{cond1}
C(0,0) = 0.
\end{equation}
Denote by $\coef_X(m)$ the coefficient at monomial $m$ in $X$. After substitutions, we have $\coef_F(z) \neq 0$ (which is true for all sextactic points and points of type $9$ on $F$). Moreover, we have $F(0,0) = 0$, so the smallest degree of a monomial in $F$ is $1$. Let
\[C_1 = C - \frac{\coef_C(z)}{\coef_F(z)}F,\]
We have $\coef_{C_1}(z) = 0$. Note that we also have
\begin{equation}\label{cond2}
\coef_{C_1}(y) = 0,
\end{equation}
because otherwise by Lemma \ref{lem: intersect} we would have $I_P(F, C) = I_P(F, C_1) \leq 1$. Therefore the smallest degree of a monomial in $C_1$ is $2$. Let
\[C_2 = C_1 - \frac{\coef_{C_1}(z^2)}{\coef_F(z)}zF, \quad C_3 = C_2 - \frac{\coef_{C_2}(yz)}{\coef_F(z)}yF.\]
We have $\coef_{C_3}(z^2) = \coef_{C_3}(yz) = 0$. We also have
\begin{equation}\label{cond3}
\coef_{C_3}(y^2) = 0,
\end{equation}
because otherwise by Lemma \ref{lem: intersect} we would have $I_P(F, C) = I_P(F, C_3) \leq 2$.

As a result, for each point on $F$, we obtain three conditions which have to be satisfied for such a curve $C$ to exist, namely \eqref{cond1}, \eqref{cond2} and \eqref{cond3}. We therefore have $6$ conditions for each pair of points. 

We verified whether such a conic exists for each pair of sextactic points and for each pair of points of type $9$ on Fermat cubic. We obtained the following results.

\begin{theorem}\label{thm: conics_s}
There are $108$ conics passing through two different sextactic points on Fermat cubic, intersecting the cubic at those points with multiplicities $3$. Moreover, for each sextactic point, there are exactly $8$ such conics passing through it.
\end{theorem}

\begin{theorem2}{\ref{thm: conics_s}'}{\label{thm: conics_t9}}
There are $324$ conics passing through two different points of type $9$ on Fermat cubic, intersecting the cubic at those points with multiplicities $3$. Moreover, for each point of type $9$, there are exactly $9$ such conics passing through it.
\end{theorem2}

Since by Proposition \ref{prop: thirdpoint} we expect $9$ conics passing through a generic point on a curve, the fact that there are only $8$ conics passing through a sextactic point may seem strange. However, this is a consequence of Theorem \ref{thm: tangent_s}: for each sextactic point $S$ on a Fermat cubic there is a flex point such that the line passing through this point and $S$ is tangent to the cubic at $S$, therefore it does not intersect the cubic at any other (sextactic) point.

For the remaining part of this section, whenever we refer to conics, we mean the conics mentioned in Theorem \ref{thm: conics_s} and Theorem \ref{thm: conics_t9}.

Two different conics $C_1$ and $C_2$, passing through a given point on $F$, intersect at exactly one other point (because their intersection index is $3$). We can potentially have $28$ and $36$ other points of intersection for a sextactic point and a point of type $9$, respectively. However, it is possible that the actual number of such points of intersection is smaller, because some of these points can have more than two conics passing through them. In fact, this is indeed the case.

\begin{theorem}\label{thm: intersections}
    The conics passing through a given sextactic point $S$ on Fermat cubic intersect in $24$ different points besides $S$, in particular:
    \begin{itemize}
        \item there are $22$ points at which $2$ such conics intersect, \vspace{-3pt}
        \item there are $2$ points at which $3$ such conics intersect, \vspace{-3pt}
        \item all $24$ points do not lie on Fermat cubic.
    \end{itemize}  
\end{theorem}

\begin{theorem2}{\ref{thm: intersections}'}
    The conics passing through a given point $T$ of type $9$ on Fermat cubic intersect in $30$ different points besides $T$, in particular:
    \begin{itemize}
        \item there are $27$ points at which $2$ such conics intersect, \vspace{-3pt}
        \item there are $3$ points at which $3$ such conics intersect, \vspace{-3pt}
        \item all $30$ points do not lie on Fermat cubic.
    \end{itemize}  
\end{theorem2}

In both cases, we gathered all obtained points of intersection and we verified how many conics overall pass through each point. We start with the points of intersection of conics tangent in sextactic points.

\begin{theorem} \label{thm: conics_inter}
    There are $540$ points not lying on Fermat cubic $F$ which are the intersection points of at least two conics passing through the same sextactic point on $F$. In particular:
    \begin{itemize}
        \item there are $486$ points at which $2$ conics intersect, \vspace{-3pt}
        \item there are $36$ points $A_1, A_2, \dots, A_{36}$ at which $6$ conics intersect, \vspace{-3pt}        
        \item there are $18$ points $B_1, B_2, \dots, B_{18}$ at which $9$ conics intersect.
    \end{itemize}  
\end{theorem}

The list of coordinates of the first $8$ points with $6$ conics passing through them is presented below. The coordinates of the remaining $28$ points can be obtained by permuting the coordinates of those $8$ points. Here, $\mu = \sqrt[3]{2}$ and $\varepsilon \in \mathbb{C}$ is a primitive root of unity of order $3$.

\vspace{-20pt}
\begin{alignat*}{2}
  A_1 &= [1 : 1 : 0], & \quad A_5 &= [1 : \varepsilon : 0], \\
  A_2 &= [1 : 1 : -\mu^2], & \quad A_6 &= [1 : \varepsilon : -\mu^2], \\
  A_3 &= [1 : 1 : -\varepsilon\mu^2], & \quad A_7 &= [1 : \varepsilon : -\varepsilon\mu^2], \\
  A_4 &= [1 : 1 : -\varepsilon^2\mu^2], & \quad A_8 &= [1 : \varepsilon : -\varepsilon^2\mu^2],  
\end{alignat*}

Similarly, we can write a list of explicit coordinates of the first $3$ points where $9$ conics intersect:
\[B_1 = [0 : 1 : -\mu], \quad B_2 = [0 : 1 : -\varepsilon\mu], \quad B_3 = [0 : 1 : -\varepsilon^2\mu].\]

It is worth noting, that points $B_1, B_2, \dots, B_{18}$ lie on three coordinate lines $x = 0$, $y = 0$ and $z = 0$ and there are $6$ such points on each line.

Moving on to points of type $9$, in this case we do not have any points of multiplicity $6$:

\begin{theorem2}{\ref{thm: conics_inter}'}
    There are $2016$ points not lying on Fermat cubic $F$ which are the intersection points of at least two conics passing through the same point of type $9$ on $F$. In particular:
    \begin{itemize}
        \item there are $1944$ points at which $2$ conics intersect, \vspace{-3pt}
        \item there are $72$ points $C_1, C_2, \dots, C_{72}$ at which $9$ conics intersect.
    \end{itemize}  
\end{theorem2}

We therefore have a new interesting set of $72$ points, lying outside Fermat cubic. Just like points of type $9$, the coordinates of these new points can also be computed and expressed in terms of $\alpha$ and $\beta$. The list of coordinates of the first $12$ points is presented below.

\vspace{-10pt}
\begin{align*}
C_1 &= [0 : 1 : \beta], & C_7 &= [0 : 3 : \alpha\beta(\beta-1)^2(\beta^2 + \beta + 1)], \\
C_2 &= [0 : 1 : \beta^2], & C_8 &= [0 : 3 : \alpha\beta(\beta^4 + 2\beta^3 + 2\beta + 1)], \\
C_3 &= [0 : 1 : \beta^4], & C_9 &= [0 : 3 : -\alpha\beta(2\beta^4 + \beta^3 + \beta + 2)], \\
C_4 &= [0 : 3 : \alpha\beta(\beta - 1)(\beta^3 + 2)], & C_{10} &= [0 : 3 : -\alpha^2\beta^2(\beta^2 + \beta + 1)], \\
C_5 &= [0 : 3 : -\alpha(\beta^2 + \beta + 1)(\beta^3 - \beta + 1)], & C_{11} &= [0 : 3 : \alpha\beta(\beta^4 - \beta^3 - \beta - 2)], \\
C_6 &= [0 : 3 : \alpha\beta(\beta^4 + 2\beta^3 - \beta + 1)], & C_{12} &= [0 : 3 : -\alpha^2\beta(\beta^2+\beta+1)(\beta^2-\beta+1)].
\end{align*}

As before, coordinates of the remaining $60$ points can be obtained by permuting the coordinates of these $12$ points above. Moreover, just like for sextactic points, the points $C_1, C_2, \dots, C_{72}$ lie on three lines $x = 0$, $y = 0$ and $z = 0$ and there are $24$ points on each line.

\section{Appendix}

The code used to obtain the results presented in this paper is available at the following link: 
\\
\url{https://github.com/maciej-zieba/fermat_cubic_special_points}.

All procedures can be run using Singular \cite{DGPS} with the following structure of the main file:

\begin{verbatim}
// ...

// Load necessary procedures for the code below
< "src/procedures.cpp";

// Compute the 9-type points on the Fermat cubic
< "src/points.cpp";

// Compute polynomial of degree 24 on which 9-type points lie
< "src/ideal_of_points.cpp";

// Compute tangents to the Fermat cubic at 6-type points
< "src/tangents_6_type.cpp";

// Compute tangents to the Fermat cubic at 9-type points
< "src/tangents_9_type.cpp";

// Compute conics passing through 6-type points with multiplicity 3
< "src/conics_6_type.cpp";

// Compute conics passing through 9-type points with multiplicity 3
< "src/conics_9_type.cpp";

\end{verbatim}

%*****************************************************************************

\bibliographystyle{abbrv}
\bibliography{main}

\end{document}